\documentclass[a4paper,11pt]{amsart}

\usepackage{amssymb}
\usepackage{epic,eepic}

\input xy

\xyoption{all}

\theoremstyle{definition}
\newtheorem{thm}{Theorem}[section]

\newtheorem{lem}[thm]{Lemma}

\newtheorem{defn}[thm]{Definition}

\newtheorem{rem}[thm]{Remark}

\newtheorem{exm}[thm]{Example}

\DeclareMathOperator{\GL}{GL}

\newcommand{\crlc}{{\mathcal C}}

\newcommand{\crls}{{\mathcal S}}
\newcommand{\crlt}{{\mathcal T}}

\newcommand{\ra}{\rightarrow}

\def\Z{\mathbb{Z}}

\def\C{\mathbb{C}}

\newcommand{\SB}[1]    {\overline{\crls} \!{#1}}
\newcommand{\vectc}{ \mathrm{Vect}_{\C}}
\newcommand{\monrep}[1]{ \mathrm{MonRep}_{\C}({#1})}
\newcommand{\frob}[1]{ \mathrm{FrobAlg}_{\C[{#1}]}}
\newcommand{\mapb}[1] {\mbox{Map}_*({#1},X)}
\newcommand{\Stilde}{\tilde{\Sigma}}

\title{Representations of the homotopy surface category of a simply
  connected space}

\begin{document}

\bibliographystyle{plain}

\vspace*{1cm}

\author{Mark Brightwell}
\address{Department of Mathematics\\
       Heriot-Watt University\\
       Edinburgh EH14 4AS\\
       Scotland}
\email{mb@ma.hw.ac.uk}
\author{Paul Turner}
\address{Department of Mathematics\\
       Heriot-Watt University\\
       Edinburgh EH14 4AS\\
       Scotland}
\email{paul@ma.hw.ac.uk}

\maketitle

\section{Introduction}
At the heart of the axiomatic formulation of 1+1-dimensional
topological field theory is the set of all surfaces with boundary 
assembled into a category. This category of surfaces has compact
1-manifolds as objects and smooth oriented cobordisms as
morphisms. Taking disjoint unions gives a monoidal
structure and a 1+1-dimensional topological field theory can be
defined to be  a monoidal 
functor $E:\crls \rightarrow \vectc$ 
normalised so that $E(\mbox{cylinder}) = \mbox{identity}$. Here
$\vectc$ is the category of finite
dimensional complex vector spaces and linear maps, with monoidal 
structure the usual tensor product. 
It is well known that 
to specify a 1+1-dimensional topological field theory
  is the same thing as specifying a finite 
dimensional commutative Frobenius algebra.

There is a natural generalisation of the category of surfaces, where
one introduces a background space $X$ and considers maps of surfaces
to $X$. Once again a monoidal category emerges. In this paper we
explain how for a simply connected background space, monoidal functors
from this category to $\vectc$ can be interpreted in terms of
Frobenius algebras with additional structure. The main result is
Theorem \ref{thm:main}. Although we use $\C$ throughout, any other
field will give the same results. 

We thank Ran Levi, Paolo Salvatore and Ulrike Tillmann for helpful
conversations.

\section{Surface categories}
There are a number of candidates for the title of surface category of a
space, all requiring some care in their formulation. Firstly recall the
definition of the category of surfaces which we denote by $\crls $
(see for example \cite{ab:tdt}, \cite{ti:csd}). 
Let $C$ be a parametrised circle. There is an object
$S_n$ of $\crls$ for each natural number $n$ consisting of an ordered set of
$n$-copies of $C$. Morphisms from $S_m$ to $S_n$ are pairs $(\Sigma,
\alpha)$ where $\Sigma$ is a smooth oriented surface equipped with an
orientation preserving homeomorphism $\alpha\colon \partial \Sigma \ra
S_m^{op} \sqcup S_n$, where $S_m^{op}$ means $S_m$ with its opposite
orientation. Two pairs are identified if there is a diffeomorphism
$T\colon \Sigma_1 \ra \Sigma_2$ satisfying $\alpha_2 \circ
T|_{\partial\Sigma} = \alpha_1$. Composition of morphisms is by gluing
surfaces via the maps $\alpha$ (in fact because of the identifications
of morphisms this is determined by a labelling of the boundary
components and their orientations). Also $\crls$ has a monoidal
structure corresponding to juxtaposition of objects and surfaces.  

It is a common procedure (e.g. in topological sigma models,
string theory etc ) to introduce a background manifold $X$ and
no longer consider abstract surfaces but maps (of a suitable kind)
from surfaces to $X$. That one can assemble a category in this way is
an idea of Segal and we refer to \cite{se:ec} for
an early appearance of this idea. Of the possible variants available
we require our surface categories to have a minimum of structure, in
particular we can deform surfaces in $X$  by homotopies relative to
the boundary. 

\begin{defn}
The
{\em homotopy surface category of } $X$, denoted by $\SB X$, is
the category defined as follows:
\begin{itemize}
\item Objects are pairs $(S,s)$ where $S$ is an object of $\crls$ 
%ordered set of $n$-copies of $C$ ($n\geq 0$) 
and $s\colon S\ra X$ is a 
continuous function.
\item Morphisms from $(S,s)$ to $(S^\prime,s^\prime)$ are triples $(\Sigma,
  \alpha, \sigma)$ where
\begin{enumerate}
\item $\Sigma$ is an smooth oriented surface
\item $\alpha\colon \partial \Sigma \ra S^{op} \sqcup S^\prime$ is an
  orientation preserving homeomorphism and
\item $\sigma\colon \Sigma \ra X$ is a continuous function such that 
$\sigma|_{\partial \Sigma}\circ \alpha^{-1} = s \sqcup s^\prime$.
\end{enumerate}
We identify $(\Sigma_1, \alpha_1, \sigma_1)$ with $(\Sigma_2,
\alpha_2, \sigma_2 )$ if there exists a diffeomorphism
$T\colon \Sigma_1 \ra \Sigma_2$ satisfying $\alpha_2 \circ
T|_{\partial\Sigma} = \alpha_1$ and $\sigma_2 \circ T \simeq_{\partial\Sigma}
\sigma_1$ (here $\simeq_{\partial\Sigma}$ means homotopy relative to
the boundary). 
\end{itemize}
\end{defn}

It is sometimes useful to think of this as the loop world equivalent
of the fundamental groupoid. Composition is by gluing using the maps
$\alpha$ for identification.
$\SB X$ has a natural monoidal structure:
 define $(S\otimes S^\prime, s\otimes s^\prime)$ and 
$(\Sigma_1 \otimes \Sigma_2, \alpha_1\otimes \alpha_2 ,
\sigma_1 \otimes \sigma_2)$ by taking 
the maps $s\otimes s^\prime$ and $\sigma_1 \otimes \sigma_2$ to be the 
disjoint union of the two original maps on the juxtaposed components.

\begin{rem} 
When $X=pt$ there is a unique map $\Sigma \ra X$ for any surface $X$ \
and we have $\SB X = \crls$.
\end{rem}

Let  $\vectc$ be the category 
of finite dimensional complex vector spaces and linear maps, with monoidal 
structure the usual tensor product and recall:

\begin{defn}
A {\em monoidal representation} of a a monoidal category $\crlc$ is a
functor $F\colon \crlc \ra \vectc$ respecting the monoidal structure
i.e. $F(0) = \C$ and $F(A\otimes B) = F(A)\otimes_{\C} F(B)$.
\end{defn}

We will write $\monrep {\SB X}$ for the set of monoidal
representations of $\SB X$. This is in fact a category where a
morphism between $F_1$ and $F_2$ is a natural transformation of
monoidal functors. It is this category for a simply connected space
which this paper investigates.

Segal suggests that just as representations of path
categories give rise to the geometry of vector bundles, so should the
representations of surface categories give rise to interesting
geometry involving the free loop space. For the  simply connected case 
and monoidal representations we will write about this elsewhere.

\section{$A$-Frobenius Algebras}

There is a ``folk'' theorem (made precise by Abrams \cite{ab:tdt}) stating
that there is a 1-1 correspondence between 1+1-dimensional topological field
theories and finite dimensional commutative Frobenius
algebras over $\C$. Recall that a Frobenius algebra over $\C$ is an
algebra $V$ over $\C$ with unit for which there exists an isomorphism of
$V$-modules $V\cong V^*$. Equivalently there exists a linear form
$\theta \colon V \ra \C$ whose kernel contains no non-trivial ideals . The
correspondence assigns the vector space $V=E(\mbox{circle})$ 
to a 1+1-dimensional topological field theory $E\colon
\crls \ra \vectc$ and the genus zero surface with 2+1 boundary circles
(the pair of pants) gives the algebra multiplication and the cap gives
the trace map $\theta$. 

 In order to generalise this to
monoidal representations of $\SB X$ for simply connected $X$ 
we need the notion of an
$A$-Frobenius algebra, where $A$ is a finitely generated abelian group. 

\begin{defn}
Let  $A$ be a finitely generated abelian group. An $A$-{\em Frobenius
  algebra} over $\C$  is a commutative $\C[A]$-algebra $V$ with unit,
  for which there exists an isomorphism of $V$-modules $V\cong V^*$,
  where $V^*=\mbox{Hom}_{\C}(V,\C)$.
\end{defn}

\begin{rem}
Alternatively one can think of an $A$-{\em Frobenius algebra} over $\C$ as a
commutative Frobenius algebra $V$ together with a representation
$i\colon A\ra \GL (V)$ such that 
\[
i(a)(x\cdot y)=i(a)x\cdot y=x\cdot i(a)y
\]
where $a\in A$, $x,y\in V$ and $x\cdot y$ is the algebra product.
\end{rem}

$A$-Frobenius algebras form a category which we will denote $\frob
A$. Morphisms are taken to be $\C[A]$-algebra isomorphisms preserving
the form $\theta$.

\begin{exm} One of the first examples of a Frobenius algebra is the
  group ring $\C[H]$ for a finite abelian group $H$ where the Frobenius
  map is given by $\theta=[1]^*$. Here $\{[h]\}_{h\in H}$ is a 
basis for the vector space $\C[H]$, and $\{[h]^*\}$ the dual basis 
of $\C[H]^*$. For details see Abrams 
\cite{ab:tdt}. Now let $A$ be a subgroup of $H$ then $\C[H]$ is an 
$A$-Frobenius algebra under the natural action of $A$ i.e. for $a\in A$
\[
i(a)(\sum_{h\in H} \alpha_h[h])=\sum_{h\in H} \alpha_h[ah]
\]
\end{exm}

\begin{exm}
Another source of Frobenius algebras is the cohomology ring $H^*(M)$
of a compact manifold where Poincar\'e duality ensures the Frobenius
structure. Quantum cohomology also provides Frobenius algebras. 
Allowing $V$ to be infinite dimensional it is known
that $QH^*(M)$ of a symplectic manifold 
is a Frobenius algebra, but it is
also naturally a $\Gamma$-Frobenius algebra, where $\Gamma$ is the
image of the Hurewicz homomorphism into integral
homology modulo torsion. Here we are taking
$QH^*(M) \cong H^*(M)\otimes \Lambda$ where $\Lambda$ is the Novikov
ring associated to the homomorphism $\omega \colon \Gamma\ra \Z$ induced
by the symplectic form (see McDuff and Salamon \cite{mcsa:jhc} for details). 
Using the notation of \cite{mcsa:jhc} an element $x\in QH^*(M)$ may be
written in the form 
\[
x= \sum_{b\in \Gamma }x_b e^{2\pi i b}
\]
and the $\Gamma$ action is given by
\[
i(a)x =  \sum_{b\in \Gamma }x_b e^{2\pi i (b+a)}
\]

\end{exm}

\section{An equivalence  of categories}

For the rest of the paper $X$ will be simply connected. The main 
Theorem of the paper is:

\begin{thm}\label{thm:main}
Let $X$ be a simply connected topological space. There is an
equivalence of categories
\[
\monrep {\SB X} \longleftrightarrow \frob {\pi_2X}
\]
\end{thm}

In order to prove this result we need to first cut down the size of
$\SB X$ by considering its skeleton, that is we identify all isomorphic
objects of the category to obtain a new equivalent category.

\begin{lem}
For $X$ simply connected, $(S,s)\cong (S^\prime, s^\prime)$ if and
only if $S=S^\prime$.
\end{lem}
\begin{proof}
Suppose that $S=S^\prime$, then since $X$ is simply connected $s\simeq
s^\prime$, that is we have a homotopy $S\times I \ra X$ agreeing with
$s$ and $s^\prime$ on the boundary. Thus taking $\Sigma=S\times I$
with the obvious identifications with $S$ at on the boundary we see we
have a morphism from  $(S,s)$ to $ (S^\prime, s^\prime)$. It is clear
this morphism is invertible and hence $(S,s)\cong (S^\prime,
s^\prime)$.

Conversely, suppose that $(S,s)\cong (S^\prime,
s^\prime)$ via a morphism $(\Sigma, \alpha, \sigma)$ with inverse
$(\Sigma^*, \alpha^*, \sigma^*)$ . Then the composition $(\Sigma,
\alpha, \sigma)(\Sigma^*, \alpha^*, \sigma^*)$ must be
identified with the identity morphism on $(S,s)$. In particular the
glued surface $\Sigma\Sigma^*$ must be diffeomorphic to $S\times I$ in
which case $\Sigma $ and $\Sigma^*$ must be made up of
cylinders and so $S=S^\prime$.
\end{proof}

We chose a representative of each isomorphism class as follows. Define
\[
c_n\colon S_n= \coprod_n S^1 \ra X
\]
to be the unique collapse map to a chosen basepoint of $X$. Define 
$\crlt X$ to be the full sub category of $\SB X$ with set of objects
$\{ (S_n,c_n)\}$. By the above lemma $\crlt X$ is the skeleton of $\SB
X$. 

Observe that 
$\crlt X$ has a natural
monoidal structure induced by that on $\SB X$. So we have:

\begin{lem}\label{lem:tx}
There is an equivalence of monoidal categories $\SB X \leftrightarrow
\crlt X$.
\end{lem}

For convenience we will denote the object $((S_n,c_n)$ simply as $n$. 
Note that using the monoidal structure $n= 1\otimes 1\otimes \cdots
\otimes 1$. 

\begin{lem}
Let $\Sigma$ be an arbitrary surface of genus $l\geq 0$ consisting of
$d$ connected components (closed or with boundary). Pick a base point
$*$ on $X$. If $X$ is simply connected then there is a 
one-to-one correspondence 
\[
[\Sigma,\partial \Sigma;X,*] \longleftrightarrow (\pi_2X)^d
\]
\end{lem}

\begin{proof}
First assume $\Sigma$ is closed, connected of genus $g$. There is a
cofibration 
\[
\xymatrix{
S^1 \ar[r]^\chi &\bigvee_{2g} S^1 \ar[r] & \Sigma
}
\]
where $\chi$ is the ``commutator'' map. Applying $\mapb -$ we get a
cofibration 
\[
\xymatrix{
\mapb \Sigma \ar[r] & (\Omega X)^{2g} \ar[r]^\chi & \Omega X
}
\]
The long exact homotopy sequence gives
\[
\xymatrix{
\cdots \ar[r] & (\pi_1\Omega X)^{2g} \ar[r]^{\chi_*} & \pi_1 \Omega X
\ar[r] & [\Sigma, X]_* \ar[r] & \pi_0(\Omega X)^{2g} \ar[r] & \cdots
}
\]
$\chi_*$ is trivial since $\Omega X\colon \Omega(\Omega X)^{2g} \ra
\Omega^2X$ is null by 
homotopy commutativity. Since $\pi_1X = \pi_0\Omega X$ is trivial it
follows that the middle map is an
equivalence. Moreover for simply connected $X$ we have $[\Sigma, X]_*
= [\Sigma, X]$ so $\pi_2X = [\Sigma, X] = [\Sigma,\partial
\Sigma;X,*]$. 

Next assume that $\Sigma$ is connected, genus $g$ and $\partial \Sigma
\neq \emptyset$ with $k$ boundary components. Let $\Stilde$ be the
surface obtained from $\Sigma$ by sewing in $k$ discs. Let $p_i$ be a
point in the $i$'th disc. There is a fibration
 \[
\xymatrix{
\mbox{Map}(\Stilde, \{p_1, \ldots, p_k\};X,*) \ar[r] & 
 \mbox{Map}(\Stilde, X) \ar[r] & X^k
}
\]
whence an exact sequence
 \[
\xymatrix{
(\pi_1X)^k \ar[r] & 
[\Stilde, \{p_1, \ldots, p_k\};X,*] \ar[r] & 
[\Stilde, X] \ar[r] & (\pi_0X)^k
}
\]
The extremities are trivial so $[\Stilde, \{p_1, \ldots, p_k\};X,*] = 
[\Stilde, X]$. Thus 
\[
[\Sigma,\partial \Sigma;X,*] = [\Stilde, \{p_1,
\ldots, p_k\};X,*] =  [\Stilde, X] = \pi_2X.
\]

Finally for a general surface each of its $d$ connected components
(closed or not) corresponds to an element of $\pi_2X$ by the above.
\end{proof}

In view of this discussion we will denote a morphism
$(\Sigma,\alpha, \sigma)$ of $\crlt X$ by $\Sigma_{\underline{g}}$ where
$\underline{g}=(g_1,\dots,g_d)\in (\pi_2X)^d$ is the 
sequence of elements of $\pi_2X$ given by the above Lemma. 
The maps
$\alpha$ will be suppressed. The unit of
$\pi_2X$ will be denoted 
by $1$ and sequence of $1$'s by $\underline{1}\in (\pi_2X)^d$. 

\begin{rem}
In other
words morphisms of $\crlt X$ are  morphisms of
$\crls$ with a labelling of each component from $\pi_2X$.
In light of this it is easy to see that if $X$ is 2-connected then
there is an equivalence of categories $\SB X \leftrightarrow \crls$.
\end{rem}

\begin{rem}
It is tempting to think there might be a monoidal equivalence between
$\crlt X$ and $\crls \times \pi_2X$. This cannot be the case in
general: any such equivalence takes $c_0$ to $(S_0,1)$ and would induce
an equivalence between the monoids of endomorphisms of these
objects, however these can be seen to be different. 
\end{rem}

\begin{lem}\label{lem1}
Let $\Sigma $ and $\Sigma^\prime $ be surfaces. Then
\[
\Sigma'_{\underline{g}} \Sigma_{\underline{h}}=(\Sigma '\Sigma)_{\underline{f}}
\]
where $f_i$ is the product of all $g_j$'s and $h_k$'s forming the connected component $f_i$ of $\Sigma' \Sigma$.
\end{lem}
\begin{proof}
Gluing of two connected components corresponds to multiplication in $\pi_2X$, and for the $i$'th component of $\Sigma' \Sigma$ we multiply all components of $\Sigma'$ and $\Sigma$ forming that component.

\vspace*{0.5cm}
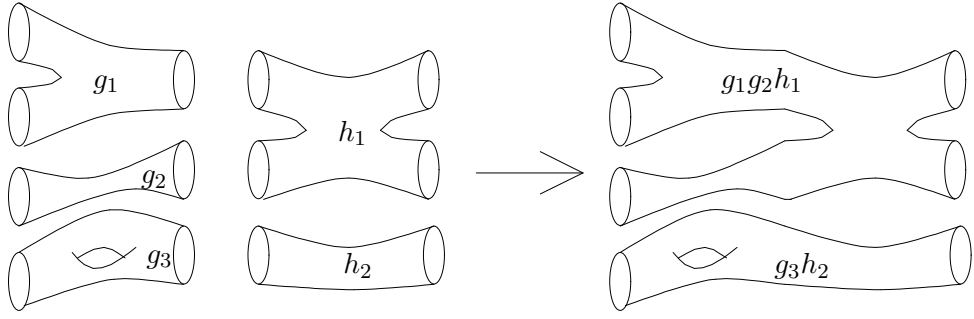
\begin{center}
\setlength{\unitlength}{0.00083333in}
\begingroup\makeatletter\ifx\SetFigFont\undefined%
\gdef\SetFigFont#1#2#3#4#5{%
  \reset@font\fontsize{#1}{#2pt}%
  \fontfamily{#3}\fontseries{#4}\fontshape{#5}%
  \selectfont}%
\fi\endgroup%
{\renewcommand{\dashlinestretch}{30}
\begin{picture}(6036,1964)(0,-10)
\put(1110,1454){\ellipse{132}{362}}
\put(1110,891){\ellipse{132}{362}}
\put(1110,360){\ellipse{132}{362}}
\put(1575,1454){\ellipse{132}{362}}
\put(1575,891){\ellipse{132}{362}}
\put(1575,360){\ellipse{132}{362}}
\put(2636,1454){\ellipse{132}{362}}
\put(2636,891){\ellipse{132}{362}}
\put(2670,360){\ellipse{132}{362}}
\put(81,1752){\ellipse{132}{362}}
\put(81,1222){\ellipse{132}{362}}
\put(81,724){\ellipse{132}{362}}
\put(81,195){\ellipse{132}{362}}
\path(81,1404)	(149.116,1410.200)
	(205.690,1416.605)
	(251.494,1423.404)
	(287.300,1430.790)
	(346.000,1470.000)

\path(346,1470)	(290.210,1518.226)
	(254.222,1530.093)
	(207.872,1542.480)
	(150.389,1555.683)
	(81.000,1570.000)

\path(115,1039)	(169.479,1062.167)
	(220.210,1083.564)
	(267.390,1103.262)
	(311.211,1121.336)
	(351.870,1137.858)
	(389.560,1152.900)
	(456.812,1178.839)
	(514.526,1199.735)
	(564.257,1216.172)
	(607.563,1228.733)
	(646.000,1238.000)

\path(646,1238)	(714.176,1249.655)
	(756.833,1254.392)
	(807.009,1258.531)
	(866.059,1262.166)
	(935.337,1265.392)
	(974.236,1266.882)
	(1016.200,1268.305)
	(1061.398,1269.674)
	(1110.000,1271.000)

\path(81,1934)	(138.290,1903.298)
	(191.680,1874.933)
	(241.377,1848.811)
	(287.589,1824.837)
	(330.524,1802.914)
	(370.390,1782.949)
	(441.748,1748.510)
	(503.324,1720.757)
	(556.782,1698.928)
	(603.787,1682.263)
	(646.000,1670.000)

\path(646,1670)	(714.435,1656.264)
	(757.156,1651.006)
	(807.354,1646.694)
	(866.382,1643.232)
	(935.596,1640.526)
	(974.446,1639.427)
	(1016.350,1638.481)
	(1061.479,1637.676)
	(1110.000,1637.000)

\path(81,907)	(125.151,896.453)
	(166.274,886.874)
	(240.069,870.523)
	(303.659,857.755)
	(358.312,848.380)
	(405.302,842.208)
	(445.899,839.047)
	(513.000,841.000)

\path(513,841)	(558.904,847.847)
	(609.520,859.666)
	(666.595,877.137)
	(731.876,900.938)
	(768.140,915.423)
	(807.110,931.746)
	(849.006,949.990)
	(894.045,970.240)
	(942.445,992.581)
	(994.425,1017.099)
	(1050.204,1043.877)
	(1110.000,1073.000)

\path(81,541)	(154.617,571.337)
	(223.254,599.167)
	(287.179,624.579)
	(346.658,647.657)
	(401.959,668.488)
	(453.350,687.159)
	(501.097,703.755)
	(545.469,718.363)
	(586.732,731.068)
	(625.154,741.959)
	(694.545,758.637)
	(755.781,769.087)
	(811.000,774.000)

\path(811,774)	(858.273,772.060)
	(919.175,761.430)
	(956.935,752.379)
	(1000.740,740.585)
	(1051.468,725.855)
	(1110.000,708.000)

\path(81,376)	(136.188,409.368)
	(187.756,440.045)
	(235.911,468.130)
	(280.862,493.721)
	(322.816,516.916)
	(361.981,537.813)
	(432.778,573.106)
	(494.914,600.386)
	(550.055,620.439)
	(599.862,634.048)
	(646.000,642.000)

\path(646,642)	(720.470,642.419)
	(764.701,636.983)
	(815.401,627.372)
	(873.927,613.278)
	(941.632,594.390)
	(979.350,583.052)
	(1019.871,570.400)
	(1063.365,556.396)
	(1110.000,541.000)

\path(115,44)	(172.745,64.500)
	(226.526,83.365)
	(276.547,100.657)
	(323.013,116.437)
	(366.131,130.764)
	(406.103,143.701)
	(477.436,165.648)
	(538.652,182.763)
	(591.391,195.535)
	(637.294,204.451)
	(678.000,210.000)

\path(678,210)	(743.346,212.925)
	(783.504,211.833)
	(830.348,208.909)
	(885.148,204.059)
	(949.176,197.188)
	(1023.703,188.200)
	(1065.301,182.882)
	(1110.000,177.000)

\path(1575,1073)	(1614.510,1074.947)
	(1650.730,1077.074)
	(1713.730,1081.961)
	(1764.867,1087.851)
	(1805.011,1094.934)
	(1873.000,1139.000)

\path(1873,1139)	(1811.116,1199.296)
	(1770.591,1215.201)
	(1718.308,1232.199)
	(1653.401,1250.672)
	(1615.941,1260.580)
	(1575.000,1271.000)

\path(2636,1073)	(2596.505,1074.947)
	(2560.298,1077.074)
	(2497.304,1081.961)
	(2446.140,1087.851)
	(2405.927,1094.934)
	(2337.000,1139.000)

\path(2337,1139)	(2399.823,1199.296)
	(2440.417,1215.201)
	(2492.725,1232.199)
	(2557.627,1250.672)
	(2595.074,1260.580)
	(2636.000,1271.000)

\path(1575,1637)	(1631.381,1614.174)
	(1683.982,1593.292)
	(1733.008,1574.290)
	(1778.665,1557.106)
	(1821.161,1541.678)
	(1860.703,1527.944)
	(1931.745,1505.305)
	(1993.444,1488.691)
	(2047.452,1477.602)
	(2095.420,1471.538)
	(2139.000,1470.000)

\path(2139,1470)	(2177.644,1472.775)
	(2220.096,1479.722)
	(2267.807,1491.342)
	(2322.227,1508.133)
	(2384.807,1530.594)
	(2456.995,1559.227)
	(2497.146,1576.013)
	(2540.243,1594.529)
	(2586.467,1614.837)
	(2636.000,1637.000)

\path(1608,708)	(1661.853,725.966)
	(1712.042,742.416)
	(1758.763,757.399)
	(1802.209,770.965)
	(1842.575,783.161)
	(1880.057,794.037)
	(1947.142,812.026)
	(2005.023,825.325)
	(2055.255,834.325)
	(2099.395,839.419)
	(2139.000,841.000)

\path(2139,841)	(2176.076,838.918)
	(2217.409,833.466)
	(2264.449,824.252)
	(2318.646,810.882)
	(2381.450,792.966)
	(2454.311,770.108)
	(2494.965,756.704)
	(2538.677,741.917)
	(2585.629,725.699)
	(2636.000,708.000)

\path(1575,541)	(1632.403,523.130)
	(1685.888,506.764)
	(1735.663,491.852)
	(1781.933,478.345)
	(1824.906,466.194)
	(1864.788,455.350)
	(1936.102,437.386)
	(1997.529,424.062)
	(2050.720,414.985)
	(2097.327,409.761)
	(2139.000,408.000)

\path(2139,408)	(2178.155,410.231)
	(2221.974,415.789)
	(2272.014,425.067)
	(2329.830,438.459)
	(2396.980,456.355)
	(2434.542,467.115)
	(2475.021,479.149)
	(2518.611,492.505)
	(2565.508,507.233)
	(2615.906,523.382)
	(2670.000,541.000)

\path(1575,177)	(1637.418,172.429)
	(1695.470,168.248)
	(1749.374,164.445)
	(1799.348,161.009)
	(1845.611,157.926)
	(1888.380,155.186)
	(1927.875,152.777)
	(1964.314,150.686)
	(2028.896,147.413)
	(2083.873,145.271)
	(2130.992,144.165)
	(2172.000,144.000)

\path(2172,144)	(2245.377,145.856)
	(2291.192,148.144)
	(2345.040,151.466)
	(2408.383,155.918)
	(2482.683,161.594)
	(2524.398,164.921)
	(2569.401,168.589)
	(2617.874,172.612)
	(2670.000,177.000)

\path(413,376)	(454.311,330.911)
	(488.052,300.800)
	(545.000,276.000)

\path(545,276)	(597.739,279.526)
	(655.470,300.671)
	(724.466,342.480)
	(765.148,372.086)
	(811.000,408.000)

\path(447,342)	(504.337,374.680)
	(548.415,395.557)
	(612.000,408.000)

\path(612,408)	(664.283,393.450)
	(699.190,373.100)
	(744.000,342.000)

\put(5921,1454){\ellipse{132}{362}}
\put(5921,891){\ellipse{132}{362}}
\put(5955,360){\ellipse{132}{362}}
\path(4859,1073)	(4898.660,1074.947)
	(4935.012,1077.074)
	(4998.239,1081.961)
	(5049.571,1087.851)
	(5089.899,1094.934)
	(5159.000,1139.000)

\path(5159,1139)	(5096.004,1199.296)
	(5055.294,1215.201)
	(5002.818,1232.199)
	(4937.683,1250.672)
	(4900.091,1260.580)
	(4859.000,1271.000)

\path(5921,1073)	(5881.479,1074.947)
	(5845.250,1077.074)
	(5782.235,1081.961)
	(5731.088,1087.851)
	(5690.940,1094.934)
	(5623.000,1139.000)

\path(5623,1139)	(5684.835,1199.296)
	(5725.364,1215.201)
	(5777.657,1232.199)
	(5842.579,1250.672)
	(5880.048,1260.580)
	(5921.000,1271.000)

\path(4859,1637)	(4915.533,1614.174)
	(4968.265,1593.292)
	(5017.406,1574.290)
	(5063.163,1557.106)
	(5105.744,1541.678)
	(5145.357,1527.944)
	(5216.511,1505.305)
	(5278.289,1488.691)
	(5332.354,1477.602)
	(5380.369,1471.538)
	(5424.000,1470.000)

\path(5424,1470)	(5462.595,1472.775)
	(5505.013,1479.722)
	(5552.704,1491.342)
	(5607.119,1508.133)
	(5669.706,1530.594)
	(5741.915,1559.227)
	(5782.081,1576.013)
	(5825.197,1594.529)
	(5871.442,1614.837)
	(5921.000,1637.000)

\path(4893,708)	(4946.826,725.966)
	(4996.993,742.416)
	(5043.694,757.399)
	(5087.124,770.965)
	(5127.479,783.161)
	(5164.951,794.037)
	(5232.030,812.026)
	(5289.918,825.325)
	(5340.171,834.325)
	(5384.346,839.419)
	(5424.000,841.000)

\path(5424,841)	(5461.026,838.918)
	(5502.324,833.466)
	(5549.343,824.252)
	(5603.534,810.882)
	(5666.345,792.966)
	(5739.226,770.108)
	(5779.897,756.704)
	(5823.628,741.917)
	(5870.602,725.699)
	(5921.000,708.000)

\path(4859,541)	(4916.552,523.130)
	(4970.169,506.764)
	(5020.058,491.852)
	(5066.427,478.345)
	(5109.485,466.194)
	(5149.438,455.350)
	(5220.865,437.386)
	(5282.371,424.062)
	(5335.620,414.985)
	(5382.275,409.761)
	(5424.000,408.000)

\path(5424,408)	(5463.106,410.231)
	(5506.890,415.789)
	(5556.908,425.067)
	(5614.718,438.459)
	(5681.875,456.355)
	(5719.445,467.115)
	(5759.936,479.149)
	(5803.543,492.505)
	(5850.459,507.233)
	(5900.880,523.382)
	(5955.000,541.000)

\path(4859,177)	(4921.568,172.429)
	(4979.749,168.248)
	(5033.762,164.445)
	(5083.826,161.009)
	(5130.158,157.926)
	(5172.978,155.186)
	(5212.503,152.777)
	(5248.951,150.686)
	(5313.493,147.413)
	(5368.351,145.271)
	(5415.271,144.165)
	(5456.000,144.000)

\path(5456,144)	(5529.872,145.856)
	(5575.845,148.144)
	(5629.806,151.466)
	(5693.229,155.918)
	(5767.587,161.594)
	(5809.327,164.921)
	(5854.353,168.589)
	(5902.849,172.612)
	(5955.000,177.000)

\put(3831,1752){\ellipse{132}{362}}
\put(3831,1222){\ellipse{132}{362}}
\put(3831,724){\ellipse{132}{362}}
\put(3831,195){\ellipse{132}{362}}
\path(2935,872)(3599,872)
\path(3599,872)(3333,741)
\path(3599,872)(3333,1006)
\path(3831,1404)	(3898.911,1410.200)
	(3955.351,1416.605)
	(4001.104,1423.404)
	(4036.953,1430.790)
	(4097.000,1470.000)

\path(4097,1470)	(4039.870,1518.226)
	(4003.839,1530.093)
	(3957.539,1542.480)
	(3900.188,1555.683)
	(3831.000,1570.000)

\path(3864,1039)	(3918.545,1062.167)
	(3969.333,1083.564)
	(4016.557,1103.262)
	(4060.409,1121.336)
	(4101.084,1137.858)
	(4138.774,1152.900)
	(4205.972,1178.839)
	(4263.549,1199.735)
	(4313.049,1216.172)
	(4356.018,1228.733)
	(4394.000,1238.000)

\path(4394,1238)	(4462.968,1249.655)
	(4505.856,1254.392)
	(4556.169,1258.531)
	(4615.273,1262.166)
	(4684.535,1265.392)
	(4723.403,1266.882)
	(4765.322,1268.305)
	(4810.464,1269.674)
	(4859.000,1271.000)

\path(3831,1934)	(3888.181,1903.298)
	(3941.473,1874.933)
	(3991.080,1848.811)
	(4037.209,1824.837)
	(4080.064,1802.914)
	(4119.849,1782.949)
	(4191.032,1748.510)
	(4252.399,1720.757)
	(4305.590,1698.928)
	(4352.243,1682.263)
	(4394.000,1670.000)

\path(4394,1670)	(4463.226,1656.264)
	(4506.179,1651.006)
	(4556.514,1646.694)
	(4615.597,1643.232)
	(4684.794,1640.526)
	(4723.613,1639.427)
	(4765.473,1638.481)
	(4810.545,1637.676)
	(4859.000,1637.000)

\path(3831,907)	(3875.042,896.453)
	(3916.068,886.874)
	(3989.705,870.523)
	(4053.170,857.755)
	(4107.722,848.380)
	(4154.622,842.208)
	(4195.129,839.047)
	(4262.000,841.000)

\path(4262,841)	(4308.025,847.847)
	(4358.728,859.666)
	(4415.855,877.137)
	(4481.154,900.938)
	(4517.413,915.423)
	(4556.371,931.746)
	(4598.244,949.990)
	(4643.253,970.240)
	(4691.614,992.581)
	(4743.547,1017.099)
	(4799.269,1043.877)
	(4859.000,1073.000)

\path(3831,541)	(3904.508,571.337)
	(3973.053,599.167)
	(4036.902,624.579)
	(4096.325,647.657)
	(4151.590,668.488)
	(4202.966,687.159)
	(4250.723,703.755)
	(4295.129,718.363)
	(4336.452,731.068)
	(4374.963,741.959)
	(4444.619,758.637)
	(4506.247,769.087)
	(4562.000,774.000)

\path(4562,774)	(4608.329,772.060)
	(4668.706,761.430)
	(4706.303,752.379)
	(4749.980,740.585)
	(4800.594,725.855)
	(4859.000,708.000)

\path(3831,376)	(3886.079,409.368)
	(3937.548,440.045)
	(3985.614,468.130)
	(4030.482,493.721)
	(4072.355,516.916)
	(4111.440,537.813)
	(4182.062,573.106)
	(4243.990,600.386)
	(4298.862,620.439)
	(4348.319,634.048)
	(4394.000,642.000)

\path(4394,642)	(4429.807,643.988)
	(4469.258,642.419)
	(4513.719,636.983)
	(4564.557,627.372)
	(4623.138,613.278)
	(4690.828,594.390)
	(4728.516,583.052)
	(4768.994,570.400)
	(4812.431,556.396)
	(4859.000,541.000)

\path(3864,44)	(3921.810,64.500)
	(3975.649,83.365)
	(4025.723,100.657)
	(4072.238,116.437)
	(4115.402,130.764)
	(4155.419,143.701)
	(4226.842,165.648)
	(4288.161,182.763)
	(4341.025,195.535)
	(4387.088,204.451)
	(4428.000,210.000)

\path(4428,210)	(4492.981,212.925)
	(4533.015,211.833)
	(4579.757,208.909)
	(4634.468,204.059)
	(4698.405,197.188)
	(4772.829,188.200)
	(4814.368,182.882)
	(4859.000,177.000)

\path(4163,376)	(4203.957,330.911)
	(4237.580,300.800)
	(4295.000,276.000)

\path(4295,276)	(4347.401,279.526)
	(4405.122,300.671)
	(4474.533,342.480)
	(4515.612,372.086)
	(4562.000,408.000)

\path(4196,342)	(4253.551,374.680)
	(4297.700,395.557)
	(4361.000,408.000)

\path(4361,408)	(4413.818,393.450)
	(4449.248,373.100)
	(4495.000,342.000)

\path(4859,708)	(4893.000,708.000)

\path(4893,708)	(4893.000,708.000)

\put(545,1404){\makebox(0,0)[lb]{\smash{{{\SetFigFont{11}{13.2}{\rmdefault}{\mddefault}{\updefault}$g_1$}}}}}
\put(845,807){\makebox(0,0)[lb]{\smash{{{\SetFigFont{11}{13.2}{\rmdefault}{\mddefault}{\updefault}$g_2$}}}}}
\put(877,309){\makebox(0,0)[lb]{\smash{{{\SetFigFont{11}{13.2}{\rmdefault}{\mddefault}{\updefault}$g_3$}}}}}
\put(2072,1073){\makebox(0,0)[lb]{\smash{{{\SetFigFont{11}{13.2}{\rmdefault}{\mddefault}{\updefault}$h_1$}}}}}
\put(2106,243){\makebox(0,0)[lb]{\smash{{{\SetFigFont{11}{13.2}{\rmdefault}{\mddefault}{\updefault}$h_2$}}}}}
\put(4461,1404){\makebox(0,0)[lb]{\smash{{{\SetFigFont{11}{13.2}{\rmdefault}{\mddefault}{\updefault}$g_1g_2h_1$}}}}}
\put(4793,243){\makebox(0,0)[lb]{\smash{{{\SetFigFont{11}{13.2}{\rmdefault}{\mddefault}{\updefault}$g_3h_2$}}}}}
\end{picture}
}

\end{center}
\vspace*{0.5cm}

\end{proof}

In particular for connected morphisms $\Sigma_{g_1}$ and  $\Sigma_{g_2}$ (where $g_i\in \pi_2X$),
\[
\Sigma'_{g_1} \Sigma_{g_2}=(\Sigma' \Sigma)_{g_1g_2}
\]

Now let
\[
E:\crlt X\rightarrow \vectc
\]
be a monoidal representation of $\crlt X$. $E(n)$ is determined by $V:=E(1)$ since
\[
E(n)=E(1\otimes \dots \otimes 1)=E(1)\otimes \dots \otimes e(1)=V^{\otimes n}
\]
To each morphism $\Sigma_{\underline{g}} \in \crlt X(m,n)$, $E$ assigns a linear map
\[
E(\Sigma_{\underline{g}}):V^{\otimes m}\rightarrow V^{\otimes n}
\]

\begin{lem}\label{lem:frobstruc}
A monoidal representation $E\colon \crlt X \ra \vectc$ induces a
Frobenius algebra structure on $E(1)$. 
\end{lem}

\begin{proof}
Let $\crlt^1 X$ be the subcategory of $\crlt X$ with objects those of 
$\crlt X$ and morphisms
\[
\{\Sigma_{\underline{1}} :\Sigma \text{ in }\crls\} \subset \ \text{morphisms of}\ \crlt X
\]
Note that the morphisms of $\crlt^1 X$ are closed under composition by
Lemma \ref{lem1}. The correspondence $\Sigma \leftrightarrow 
\Sigma_{\underline{1}}$ then induces an equivalence of categories
$\crls \leftrightarrow \crlt^1 X$. The composition $\crls
\leftrightarrow \crlt^1 X \ra \crlt X \ra \vectc$ gives a Frobenius
algebra structure on $E(1)$ by Proposition 13 in \cite{ab:tdt}.
\end{proof}

%Note that the algebra structure above is induced by $Y_1 \in \crlt X
%(2,1)$ where $Y\in \crls (2,1)$ is a pair-or-pants surface.

The remainder of this section is devoted to the proof of Theorem
\ref{thm:main}. 

First we see that a monoidal representation $E\colon \crlt X \ra
\vectc$ gives rise to a $\pi_2X$-Frobenius algebra. Let $V=E(1)$ which
courtesy of Lemma \ref{lem:frobstruc} is a Frobenius algebra. Denote
by $I_g\in \crlt X(1,1)$ the cylinder labelled with $g\in
\pi_2X$. Define 
\[
i\colon \pi_2X \ra GL(V)
\]
by 
\[
i(g) = E(I_g)
\]
Writing $i_g$ for $i(g)$, note that $i_g$ is indeed invertible since
\[
E(I_g)E(I_{g^{-1}})=E(I_gI_{g^{-1}})=E(I_1)=\mathrm{id}_V.
\]
The property 
\begin{equation}\label{eq:condition}
i_g(x\cdot y)=i_gx\cdot y=x\cdot i_gy.
\end{equation}
follows from

\vspace*{0.5cm}
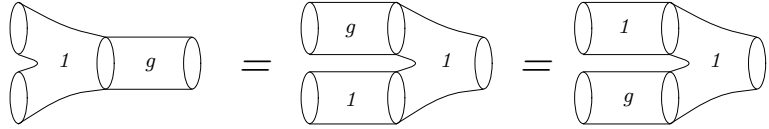
\begin{center}
\setlength{\unitlength}{0.00050000in}
\begingroup\makeatletter\ifx\SetFigFont\undefined%
\gdef\SetFigFont#1#2#3#4#5{%
  \reset@font\fontsize{#1}{#2pt}%
  \fontfamily{#3}\fontseries{#4}\fontshape{#5}%
  \selectfont}%
\fi\endgroup%
{\renewcommand{\dashlinestretch}{30}
\begin{picture}(7888,1310)(0,-10)
\put(105,1009){\ellipse{180}{542}}
\put(105,286){\ellipse{180}{542}}
\put(1911,647){\ellipse{180}{542}}
\put(1008,647){\ellipse{180}{542}}
\put(3131,1009){\ellipse{180}{542}}
\put(4938,647){\ellipse{180}{542}}
\put(3131,286){\ellipse{180}{542}}
\put(4034,286){\ellipse{180}{542}}
\put(4034,1009){\ellipse{180}{542}}
\put(5977,1009){\ellipse{180}{542}}
\put(7783,647){\ellipse{180}{542}}
\put(5977,286){\ellipse{180}{542}}
\put(6880,286){\ellipse{180}{542}}
\put(6880,1009){\ellipse{180}{542}}
\spline(105,557)
(376,647)(105,738)
\spline(105,15)
(601,331)(1008,376)
\spline(105,1280)
(601,1009)(1008,918)
\path(1008,918)(1911,918)
\path(1008,376)(1911,376)
\path(3131,1280)(4034,1280)
\path(3131,738)(4034,738)
\spline(4034,557)
(4305,647)(4034,738)
\spline(4034,15)
(4531,331)(4938,376)
\spline(4034,1280)
(4531,1009)(4938,918)
\path(3131,557)(4034,557)
\path(3131,15)(4034,15)
\path(5977,1280)(6880,1280)
\path(5977,738)(6880,738)
\spline(6880,557)
(7151,647)(6880,738)
\spline(6880,15)
(7377,331)(7783,376)
\spline(6880,1280)
(7377,1009)(7783,918)
\path(5977,557)(6880,557)
\path(5977,15)(6880,15)
\put(6338,964){\makebox(0,0)[lb]{\smash{{{\SetFigFont{8}{9.6}{\rmdefault}{\mddefault}{\itdefault}1}}}}}
\put(6383,241){\makebox(0,0)[lb]{\smash{{{\SetFigFont{8}{9.6}{\rmdefault}{\mddefault}{\itdefault}g}}}}}
\put(7287,602){\makebox(0,0)[lb]{\smash{{{\SetFigFont{8}{9.6}{\rmdefault}{\mddefault}{\itdefault}1}}}}}
\put(3492,964){\makebox(0,0)[lb]{\smash{{{\SetFigFont{8}{9.6}{\rmdefault}{\mddefault}{\itdefault}g}}}}}
\put(3492,196){\makebox(0,0)[lb]{\smash{{{\SetFigFont{8}{9.6}{\rmdefault}{\mddefault}{\itdefault}1}}}}}
\put(4486,602){\makebox(0,0)[lb]{\smash{{{\SetFigFont{8}{9.6}{\rmdefault}{\mddefault}{\itdefault}1}}}}}
\put(511,602){\makebox(0,0)[lb]{\smash{{{\SetFigFont{8}{9.6}{\rmdefault}{\mddefault}{\itdefault}1}}}}}
\put(1414,602){\makebox(0,0)[lb]{\smash{{{\SetFigFont{8}{9.6}{\rmdefault}{\mddefault}{\itdefault}g}}}}}
\put(2369,500){\makebox(0,0)[lb]{\smash{{{\SetFigFont{17}{20.4}{\rmdefault}{\mddefault}{\itdefault}=}}}}}
\put(5294,500){\makebox(0,0)[lb]{\smash{{{\SetFigFont{17}{20.4}{\rmdefault}{\mddefault}{\itdefault}=}}}}}
\end{picture}
}

\end{center}
\vspace*{0.5cm}

\noindent
recalling that the pair-of-pants induces the algebra structure on
$V$. 

Moreover a natural transformation of monoidal functors induces a
morphism of $\pi_2X$-Frobenius algebras. Consider a monoidal natural
transformation $\phi$ between $E_1$ and $E_2$ which by definition
assigns to object $n$ a map $\phi_n\colon E_1(1)^{\otimes n} \ra
E_2(1)^{\otimes n}$ such that every morphism $\Sigma_{\underline{g}}$
gives a commutative diagram
\[
\xymatrix{ E_1(1)^{\otimes m} \ar[r]^-{\phi_m}
  \ar[d]_-{E_1(\Sigma_{\underline{g}})}& E_2(1)^{\otimes m}
\ar[d]^-{E_2(\Sigma_{\underline{g}})}\\
E_1(1)^{\otimes n} \ar[r]_-{\phi_n}&  E_2(1)^{\otimes n}} 
\]
By the definition of a monoidal functor we have
$\phi_n=\phi_1^{\otimes n}$ and it
follows from \cite{ab:tdt} that $\phi_1$ is a Frobenius algebra
isomorphism. $\pi_2X$-linearity follows from the commutativity of
\[
\xymatrix{ E_1(1) \ar[r]^-{\phi_1}
  \ar[d]_-{E_1(I_g)}& E_2(1)
\ar[d]^-{E_2(I_g)}\\
E_1(1) \ar[r]_-{\phi_1}&  E_2(1)} 
\]

We have thus defined a functor
\[
F\colon \monrep {\crlt X} \ra \frob {\pi_2X}.
\]

Now consider a Frobenius algebra $V$ together with
a representation $i$ of $\pi_2X$ on $V$ satisfying
(\ref{eq:condition}).  To such a Frobenius algebra we wish to assign a
monoidal representation $E$ of $\crlt X$. We will set $E(1) = V$ and 
define $E(\mbox{pair-of-pants})\colon V\otimes V \ra V$ to be the multiplication in
$V$. Letting $I_g^m=I_{g_1}\otimes \dots \otimes I_{g_m}$ be a
morphism consisting of $m$ straight cylinders such that $\prod g_i=g$
($g,g_i\in \pi_2X$), define $E(I_g^m)= i_{g_1}\otimes \cdots \otimes
i_{g_m}$, recalling that $i_{g_k} \in GL(V)$. 
To define $E$ for an arbitrary surface we
first show that  the entire category $\crlt X$ can be recovered from $\crlt^1
X$ and $\{I_h\}_{h\in \pi_2X}$ (by composition and monoidal
structure). 

First assume $\Sigma_g\in \crlt X(m,n)$ is a connected
surface and not a morphism to the empty manifold (that is
$n>0$).  Then the following diagram commutes 
\[
\xymatrix{m\ar[r]^-{\Sigma_g} \ar[dr]_-{\Sigma_1} &n \\ & n \ar[u]_-{I_g^n}}
\]
If $m>0$ we have similarly
\[
\xymatrix{m\ar[r]^-{\Sigma_g} \ar[d]_-{I_g^m} &n \\ m \ar[ur]_-{\Sigma_1}}
\]
For a closed connected surface $\Sigma_g$ $(n=m=0)$, capping the
cylinder $I_g^1$ suffices; let $D,E$ be the discs to and from the
empty manifold respectively and $\overline{\Sigma_g}$ be the morphism
in $\crlt X(1,1)$ obtained from $\Sigma_g$ by removing a disc at each
end. We have 
\[
\xymatrix{0 \ar[r]^-{\Sigma_g} \ar[d]_-{I_gE_1}&0\\1
  \ar[r]_-{\overline{\Sigma_1}}& 1 \ar[u]_-{D_1}} 
\]
It follows from this discussion that any surface
 $\Sigma_{\underline{g}}\in
\crlt X(m,n)$ can be decomposed (non-uniquely) 
into surfaces in $\crlt^1 X$ and
$\{I^m_g\}_{g\in \pi_2X}$. Using such a decomposition we can define a
 map
\[
E(\Sigma_{\underline{g}})\colon V^{\otimes m} \ra V^{\otimes n}
\]
A priori this depends on the decomposition, but as we now see it is
in fact independent of that choice. The hard work is contained in the
proof  that $\monrep{\crls} \leftrightarrow 
\mathrm{FrobAlg}_{\C}$ (see \cite{ab:tdt})
where it is shown that any decomposition of a surface into the five
basic surfaces is independent of the slicing. We need only in addition
to know that the map $E(\Sigma_{\underline{g}})$ above is independent
of the labellings of the 
components by $\pi_2X$. This follows from the fact that 
$i_g(x\cdot y)=i_gx\cdot y=x\cdot i_gy$.

Furthermore, a morphism $\phi$ of $\pi_2X$-Frobenius algebras induces a
monoidal natural transformation by setting $\phi_n=\phi^{\otimes
  n}$. The discussion above and $\pi_2X$-linearity ensure we get a
monoidal natural transformation.

So we have a functor
\[
G\colon \frob {\pi_2X}\ra \monrep {\crlt X}
\]
It is easy to see that $F\circ G =1$ and $G\circ F=1$ so there is an
isomorphism of categories
\[
\monrep {\crlt X} \cong \frob {\pi_2X}
\]
Finally using Lemma \ref{lem:tx} we get the desired equivalence
\[
\monrep {\SB X} \leftrightarrow \frob {\pi_2X}
\]
completing the proof of the main theorem.

%\bibliography{/home/paul/bibtex/topology}

\end{document}